\documentclass[a4paper,12pt]{article}
\usepackage{amsmath}
\usepackage{ifthen,amssymb,theorem,color,enumerate}
\usepackage{geometry,latexsym}
\usepackage{citesort}
\ifx\fiverm\undefined\newfont\fiverm{cmr5}\fi
\input{m-pictex}
\newdimen\einheit
\newbox\picdiagram
\def\pictexdiagram #1/#2 #3/#4 #5 {%
\def\empty{}%
\def\test{#3}\ifx\test\empty\def\hor{0mm}\else\def\hor{#3mm}\fi%
\def\test{#4}\ifx\test\empty\def\ver{0mm}\else\def\ver{#4mm}\fi%
\noindent\begin{center}%
\mbox{#1\phantom{#2}}\hfill%
\setbox\picdiagram=\hbox{\vbox{\vspace*{\ver}%
\beginpicture%
\setcoordinatesystem units <\einheit,\einheit> point at 0 0 %
#5%
\endpicture%
\def\test{#4}\ifx\test\empty\def\ver{0mm}\else\def\ver{-#4mm}\fi%
\vspace*{\ver}}}%
\parbox{\wd\picdiagram}{\makebox[\wd\picdiagram][l]{\hspace*{\hor}\box\picdiagram}}%
\hfill\mbox{\phantom{#1}#2}%
\end{center}%
}%
\def\hdeg{\plot 0 0  0.5 0.5  1 0  1.5 -0.5  2 0 /}
\def\vdeg{\plot 0 0  0.5 0.5  0 1  -0.5 1.5  0 2 /}
\def\backsquarelabels #1&#2&#3&#4 {\def\backl{#1}\def\backt{#2}\def\backr{#3}\def\backb{#4}}%
\def\middlearrowlabels #1&#2&#3&#4 {\def\middletl{#1}\def\middletr{#2}\def\middlebr{#3}\def\middlebl{#4}}%
\def\frontsquarelabels #1&#2&#3&#4 {\def\frontl{#1}\def\frontt{#2}\def\frontr{#3}\def\frontb{#4}}%
\def\facelabels #1&#2&#3&#4&#5&#6 {\def\backface{#1}\def\frontface{#2}\def\leftface{#3}\def\rightface{#4}%
\def\bottomface{#5}\def\topface{#6}}%
\def\letteredcube #1&#2&#3&#4 {%
\einheit=\zoomfactor mm%
\pictexdiagram #1/#2 #3/#4 {%
\fontsize{\fontgroesse}{\zeilenabstand}\selectfont%
\ifnum\showindexing=0\relax\else\arrow <2\einheit> [0.25,1] from 110 65 to 110 55 %
\arrow <2\einheit> [0.25,1] from 110 65 to 120 65 %
\arrow <2\einheit> [0.25,1] from 110 65 to 120 60.83 %
\put {\footnotesize $1$} [t] at 110 53 \put {\footnotesize $3$} [l] at 122 65 \put {\footnotesize $2$} [l] at 122 60 \fi%
\setplotsymbol ({.})%
\put {$\backl$} [r] at -2.5 45 \arrow <2.4\einheit> [0.25,1] from 0 43.6 to 0 43.5 \plot 0 65  0 25 /
\put {$\frontl$} [l] at 62.5 20 \arrow <2.4\einheit> [0.25,1] from 60 18.6 to 60 18.5 \plot 60 40  60 0 /
\put {$\middlebl$} [t] at 30 10 \arrow <2.4\einheit> [0.25,1] from 31.5 11.8417 to 31.6 11.8 \plot 0 25  60 0 /
\put {$\bullet$} at 0 25 \put {$$} [t] at 0 22 
\put {$\bullet$} at 60 0 \put {$$} [t] at 60 -3 
\put {$\backb$} [t] at 20 23 \arrow <2.4\einheit> [0.25,1] from 21.4 25 to 21.5 25 %
\setdots <1.4\einheit> \plot 0 25  40 25 /\setsolid
\put {$\frontb$} [t] at 80 -2 \arrow <2.4\einheit> [0.25,1] from 81.4 0 to 81.5 0 \plot 60 0  100 0 /
\put {$\middlebr$} [t] at 70 10 \arrow <2.4\einheit> [0.25,1] from 71.5 11.8417 to 71.6 11.8 %
\setdots <1.4\einheit> \plot 40 25  100 0 /\setsolid
\put {$\bullet$} at 40 25 \put {$$} [t] at 40 22 
\put {$\bullet$} at 100 0 \put {$$} [t] at 100 -3 
\put {$\backr$} [r] at 37.5 45 \arrow <2.4\einheit> [0.25,1] from 40 43.6 to 40 43.5 %
\setdots <1.4\einheit> \plot 40 65  40 25 /\setsolid
\put {$\frontr$} [l] at 102.5 20 \arrow <2.4\einheit> [0.25,1] from 100 18.6 to 100 18.5 \plot 100 40  100 0 /
\put {$\middletr$} [b] at 70 55 \arrow <2.4\einheit> [0.25,1] from 71.5 51.8417 to 71.6 51.8 \plot 40 65  100 40 /
\put {$\bullet$} at 40 65 \put {$$} [b] at 40 68 
\put {$\bullet$} at 100 40 \put {$$} [b] at 100 43 
\put {$\backt$} [b] at 20 67 \arrow <2.4\einheit> [0.25,1] from 21.4 65 to 21.5 65 \plot 0 65  40 65 /
\put {$\frontt$} [b] at 80 42 \arrow <2.4\einheit> [0.25,1] from 81.4 40 to 81.5 40 \plot 60 40  100 40 /
\put {$\bullet$} at 0 65 \put {$$} [b] at 0 68 
\put {$\bullet$} at 60 40 \put {$$} [b] at 60 43 
\put {$\middletl$} [b] at 30 55 \arrow <2.4\einheit> [0.25,1] from 31.5 51.8417 to 31.6 51.8 \plot 0 65  60 40 /
\setplotsymbol ({\fiverm .})\setdashes <1.4\einheit> %
\plot -16 32.5  30 32.5 /\put {$\leftface$} [r] at -18.5 32.5 
\plot 80 20  115 5.417 /\put {$\frontface$} [l] at 117.5 4.17 
\plot 50 52.5  50 80 /\put {$\topface$} [b] at 50 81 
\plot 50 -15  50 12.5 /\put {$\bottomface$} [t] at 50 -17 
\plot -15 59.583  20 45 /\put {$\backface$} [r] at -17.5 60.83 
\plot 116 32.5  70 32.5 /\put {$\rightface$} [l] at 118.5 32.5 
} %
}%
\def\verticalidentity#1{%
\def\test{#1}\ifx\test\empty\else\def\zoomfactor{#1}\einheit=\zoomfactor mm\fi%
\beginpicture%
\setcoordinatesystem units <\einheit,\einheit> point at 0 0 %
\linethickness1pt\linethickness=\zoomfactor\linethickness%
\putrule from 0 0 to 0 5 \putrule from 5 0 to 5 5 %
\endpicture%
}%
\def\horizontalidentity#1{%
\def\test{#1}\ifx\test\empty\else\def\zoomfactor{#1}\einheit=\zoomfactor mm\fi%
\beginpicture%
\setcoordinatesystem units <\einheit,\einheit> point at 0 0 %
\linethickness1pt\linethickness=\zoomfactor\linethickness%
\putrule from 0 0 to 5 0 \putrule from 0 5 to 5 5 %
\endpicture%
}%
\def\gammaplus#1{%
\def\test{#1}\ifx\test\empty\else\def\zoomfactor{#1}\einheit=\zoomfactor mm\fi%
\beginpicture%
\setcoordinatesystem units <\einheit,\einheit> point at 0 0 %
\linethickness1pt\linethickness=\zoomfactor\linethickness%
\putrule from 0 0 to 0 5 \putrule from 0 5 to 5 5 %
\endpicture%
}%
\def\gammaminus#1{%
\def\test{#1}\ifx\test\empty\else\def\zoomfactor{#1}\einheit=\zoomfactor mm\fi%
\beginpicture%
\setcoordinatesystem units <\einheit,\einheit> point at 0 0 %
\linethickness1pt\linethickness=\zoomfactor\linethickness%
\putrule from 0 0 to 5 0 \putrule from 5 0 to 5 5 %
\endpicture%
}%

\def\twocubeleft #1&#2&#3 {\def\ltwodownleft{#1}\def\ltwodownmiddle{#2}\def\ltwoupmiddle{#3}}%
\def\onecubeleft #1&#2&#3&#4&#5&#6 {\def\loneleft{#1}\def\loneupmiddle{#2}\def\loneupright{#3}%
\def\loneright{#4}\def\lonedownmiddle{#5}\def\lonedownleft{#6}}%
\def\onecuberight #1&#2&#3&#4&#5&#6 {\def\roneleft{#1}\def\roneupmiddle{#2}\def\roneupright{#3}\def\roneright{#4}%
\def\ronedownmiddle{#5}\def\ronedownleft{#6}}%
\def\twocuberight #1&#2&#3 {\def\rtwodownmiddle{#1}\def\rtwoupmiddle{#2}\def\rtwoupright{#3}}%
\def\letteredrectangles #1&#2&#3&#4&#5 {%
\einheit=\zoomfactor mm%
\pictexdiagram #1/#2 #3/#4 {%
\fontsize{\fontgroesse}{\zeilenabstand}\selectfont%
\ifnum\orientation=0\ifnum\showindexing=0\relax\else%
\arrow <2\einheit> [0.25,1] from 150 40 to 150 30 \arrow <2\einheit> [0.25,1] from 150 40 to 160 40 %
\put {$1$} [t] at 150 28 \put {$2$} [l] at 162 40 \fi%
\arrow <2\einheit> [0.25,1] from 0 28.8 to 0 28.7 \arrow <2\einheit> [0.25,1] from 0 8.8 to 0 8.7 \plot 0 0  0 40 /%
\setquadratic\multiput {\vdeg} at 0 21 *9 0 2 /\setlinear%
\arrow <2\einheit> [0.25,1] from 20 28.8 to 20 28.7 \arrow <2\einheit> [0.25,1] from 20 8.8 to 20 8.7 \plot 20 0  20 40 /%
\arrow <2\einheit> [0.25,1] from 40 28.8 to 40 28.7 \arrow <2\einheit> [0.25,1] from 40 8.8 to 40 8.7 \plot 40 0 40 40 /%
\arrow <2\einheit> [0.25,1] from 60 28.8 to 60 28.7 \arrow <2\einheit> [0.25,1] from 60 8.8 to 60 8.7 \plot 60 0  60 40 /%
\setquadratic\multiput {\vdeg} at 60 21 *9 0 2 /\setlinear%
\arrow <2\einheit> [0.25,1] from 11.2 0 to 11.3 0 \arrow <2\einheit> [0.25,1] from 31.2 0 to 31.3 0 %
\arrow <2\einheit> [0.25,1] from 51.2 0 to 51.3 0 \plot 0 0  60 0 /\setquadratic\multiput {\hdeg} at 40 0 *9 2 0 /%
\setlinear%
\arrow <2\einheit> [0.25,1] from 11.2 20 to 11.3 20 \arrow <2\einheit> [0.25,1] from 31.2 20 to 31.3 20 %
\arrow <2\einheit> [0.25,1] from 51.2 20 to 51.3 20 \plot 0 20  60 20 /\setquadratic\multiput {\hdeg} at 40 20 *9 2 0 /%
\setlinear%
\arrow <2\einheit> [0.25,1] from 11.2 40 to 11.3 40 \arrow <2\einheit> [0.25,1] from 31.2 40 to 31.3 40 %
\arrow <2\einheit> [0.25,1] from 51.2 40 to 51.3 40 \plot 0 40  60 40 /\setquadratic\multiput {\hdeg} at 0 40 *9 2 0 /%
\setlinear%
\linethickness1pt\linethickness=\zoomfactor\linethickness%
\put {$\gammaplus{}$} at 10 30 \put {$\gammaminus{}$} at 50 30 \put {$\horizontalidentity{}$} at 50 10 %
\put {$\ltwodownleft$} at 10 10 \put {$\ltwodownmiddle$} at 30 10 \put {$\ltwoupmiddle$} at 30 30 %
\put {$\loneleft$} [r] at -2.5 10 \put {$\lonedownleft$} [t] at 10 -2 \put {$\lonedownmiddle$} [t] at 30 -2 %
\put {$\loneupmiddle$} [b] at 30 42 \put {$\loneupright$} [b] at 50 42 \put {$\loneright$} [l] at 62.5 10 %
\put {$#5$} at 70 20 %
\arrow <2\einheit> [0.25,1] from 80 28.8 to 80 28.7 \arrow <2\einheit> [0.25,1] from 80 8.8 to 80 8.7 \plot 80 0  80 40 /%
\setquadratic\multiput {\vdeg} at 80 1 *9 0 2 /\setlinear%
\arrow <2\einheit> [0.25,1] from 100 28.8 to 100 28.7 \arrow <2\einheit> [0.25,1] from 100 8.8 to 100 8.7 \plot 100 0  100 40 /%
\arrow <2\einheit> [0.25,1] from 120 28.8 to 120 28.7 \arrow <2\einheit> [0.25,1] from 120 8.8 to 120 8.7 \plot 120 0  120 40 /%
\arrow <2\einheit> [0.25,1] from 140 28.8 to 140 28.7 \arrow <2\einheit> [0.25,1] from 140 8.8 to 140 8.7 \plot 140 0  140 40 /%
\setquadratic\multiput {\vdeg} at 140 1 *9 0 2 /\setlinear%
\arrow <2\einheit> [0.25,1] from 91.2 0 to 91.3 0 \arrow <2\einheit> [0.25,1] from 111.2 0 to 111.3 0 %
\arrow <2\einheit> [0.25,1] from 131.2 0 to 131.3 0 \plot 80 0  140 0 /\setquadratic\multiput {\hdeg} at 120 0 *9 2 0 /%
\setlinear%
\arrow <2\einheit> [0.25,1] from 91.2 20 to 91.3 20 \arrow <2\einheit> [0.25,1] from 111.2 20 to 111.3 20 %
\arrow <2\einheit> [0.25,1] from 131.2 20 to 131.3 20 \plot 80 20  140 20 /\setquadratic\multiput {\hdeg} at 80 20 *9 2 0 /%
\setlinear%
\arrow <2\einheit> [0.25,1] from 91.2 40 to 91.3 40 \arrow <2\einheit> [0.25,1] from 111.2 40 to 111.3 40 %
\arrow <2\einheit> [0.25,1] from 131.2 40 to 131.3 40 \plot 80 40  140 40 /\setquadratic\multiput {\hdeg} at 80 40 *9 2 0 /%
\setlinear%
\linethickness1pt\linethickness=\zoomfactor\linethickness%
\put {$\gammaplus{}$} at 90 10 \put {$\horizontalidentity{}$} at 90 30 \put {$\gammaminus{}$} at 130 10 %
\put {$\rtwodownmiddle$} at 110 10 \put {$\rtwoupmiddle$} at 110 30 \put {$\rtwoupright$} at 130 30 %
\put {$\roneleft$} [r] at 77.5 30 \put {$\ronedownleft$} [t] at 90 -2 \put {$\ronedownmiddle$} [t] at 110 -2 %
\put {$\roneupmiddle$} [b] at 110 42 \put {$\roneupright$} [b] at 130 42 \put {$\roneright$} [l] at 142.5 30 %
\else\ifnum\showindexing=0\relax\else%
\arrow <2\einheit> [0.25,1] from 120 40 to 120 30 \arrow <2\einheit> [0.25,1] from 120 40 to 130 40 %
\put {$1$} [t] at 120 28 \put {$2$} [l] at 132 40 \fi%
\arrow <2\einheit> [0.25,1] from 11.2 0 to 11.3 0 \arrow <2\einheit> [0.25,1] from 31.2 0 to 31.3 0 \plot 0 0  40 0 /%
\setquadratic\multiput {\hdeg} at 0 0 *9 2 0 /\setlinear%
\arrow <2\einheit> [0.25,1] from 11.2 20 to 11.3 20 \arrow <2\einheit> [0.25,1] from 31.2 20 to 31.3 20 \plot 0 20  40 20 /%
\arrow <2\einheit> [0.25,1] from 11.2 40 to 11.3 40 \arrow <2\einheit> [0.25,1] from 31.2 40 to 31.3 40 \plot 0 40  40 40 /%
\arrow <2\einheit> [0.25,1] from 11.2 60 to 11.3 60 \arrow <2\einheit> [0.25,1] from 31.2 60 to 31.3 60 \plot 0 60  40 60 /%
\setquadratic\multiput {\hdeg} at 0 60 *9 2 0 /\setlinear%
\arrow <2\einheit> [0.25,1] from 0 8.8 to 0 8.7 \arrow <2\einheit> [0.25,1] from 0 28.8 to 0 28.7 %
\arrow <2\einheit> [0.25,1] from 0 48.8 to 0 48.7 \plot 0 0  0 60 /%
\setquadratic\multiput {\vdeg} at 0 41 *9 0 2 /\setlinear%
\arrow <2\einheit> [0.25,1] from 20 8.8 to 20 8.7 \arrow <2\einheit> [0.25,1] from 20 28.8 to 20 28.7 %
\arrow <2\einheit> [0.25,1] from 20 48.8 to 20 48.7 \plot 20 0  20 60 /%
\setquadratic\multiput {\vdeg} at 20 1 *9 0 2 /\setlinear%
\arrow <2\einheit> [0.25,1] from 40 8.8 to 40 8.7 \arrow <2\einheit> [0.25,1] from 40 28.8 to 40 28.7 %
\arrow <2\einheit> [0.25,1] from 40 48.8 to 40 48.7 \plot 40 0  40 60 /%
\setquadratic\multiput {\vdeg} at 40 1 *9 0 2 /\setlinear%
\linethickness1pt\linethickness=\zoomfactor\linethickness%
\put {$\gammaminus{}$} at 10 10 \put {$\gammaplus{}$} at 10 50 \put {$\verticalidentity{}$} at 30 10 %
\put {$\ltwodownleft$} at 10 30 \put {$\ltwodownmiddle$} at 30 30 \put {$\ltwoupmiddle$} at 30 50 %
\put {$\lonedownleft$} [r] at -2.5 10 \put {$\loneleft$} [r] at -2.5 30 \put {$\lonedownmiddle$} [t] at 30 -2 %
\put {$\loneupmiddle$} [b] at 30 62 \put {$\loneright$} [l] at 42 30 \put {$\loneupright$} [l] at 42 50 %
\put {$=$} at 50 30 %
\arrow <2\einheit> [0.25,1] from 71.2 0 to 71.3 0 \arrow <2\einheit> [0.25,1] from 91.2 0 to 91.3 0 \plot 60 0  100 0 /%
\setquadratic\multiput {\hdeg} at 80 0 *9 2 0 /\setlinear%
\arrow <2\einheit> [0.25,1] from 71.2 20 to 71.3 20 \arrow <2\einheit> [0.25,1] from 91.2 20 to 91.3 20 \plot 60 20  100 20 /%
\arrow <2\einheit> [0.25,1] from 71.2 40 to 71.3 40 \arrow <2\einheit> [0.25,1] from 91.2 40 to 91.3 40 \plot 60 40  100 40 /%
\arrow <2\einheit> [0.25,1] from 71.2 60 to 71.3 60 \arrow <2\einheit> [0.25,1] from 91.2 60 to 91.3 60 \plot 60 60  100 60 /%
\setquadratic\multiput {\hdeg} at 80 60 *9 2 0 /\setlinear%
\arrow <2\einheit> [0.25,1] from 60 8.8 to 60 8.7 \arrow <2\einheit> [0.25,1] from 60 28.8 to 60 28.7 %
\arrow <2\einheit> [0.25,1] from 60 48.8 to 60 48.7 \plot 60 0  60 60 /%
\setquadratic\multiput {\vdeg} at 60 41 *9 0 2 /\setlinear%
\arrow <2\einheit> [0.25,1] from 80 8.8 to 80 8.7 \arrow <2\einheit> [0.25,1] from 80 28.8 to 80 28.7 %
\arrow <2\einheit> [0.25,1] from 80 48.8 to 80 48.7 \plot 80 0  80 60 /%
\setquadratic\multiput {\vdeg} at 80 41 *9 0 2 /\setlinear%
\arrow <2\einheit> [0.25,1] from 100 8.8 to 100 8.7 \arrow <2\einheit> [0.25,1] from 100 28.8 to 100 28.7 %
\arrow <2\einheit> [0.25,1] from 100 48.8 to 100 48.7 \plot 100 0  100 60 /%
\setquadratic\multiput {\vdeg} at 100 1 *9 0 2 /\setlinear%
\put {$\verticalidentity{}$} at 70 50 \put {$\gammaminus{}$} at 90 10 \put {$\gammaplus{}$} at 90 50 %
\put {$\rtwodownmiddle$} at 70 10 \put {$\rtwoupmiddle$} at 70 30 \put {$\rtwoupright$} at 90 30 %
\put {$\ronedownleft$} [r] at 57.5 10 \put {$\roneleft$} [r] at 57.5 30 \put {$\ronedownmiddle$} [t] at 70 -2 %
\put {$\roneupmiddle$} [b] at 70 62 \put {$\roneright$} [l] at 102 30 \put {$\roneupright$} [l] at 102 50 %
\fi%
} %
}%
\def\squarelabels #1&#2&#3&#4&#5&#6&#7&#8&#9 {\def\zerotl{#1}\def\zerotr{#2}\def\zerobr{#3}\def\zerobl{#4}%
\def\onel{#5}\def\onet{#6}\def\oner{#7}\def\oneb{#8}\def\twoc{#9}}%
\def\letteredsquare #1&#2&#3&#4 {%
\einheit=\zoomfactor mm%
\pictexdiagram #1/#2 #3/#4 {%
\fontsize{\fontgroesse}{\zeilenabstand}\selectfont%
\ifnum\showindexing=0\relax\else%
\arrow <1.25\einheit> [0.25,1] from 30 20 to 30 15 \arrow <1.25\einheit> [0.25,1] from 30 20 to 35 20 %
\put {\footnotesize $1$} [t] at 30 13.5 \put {\footnotesize $2$} [l] at 36.5 20 \fi%
\put {\scriptsize $\bullet$} at -10 0 \put {$\zerobl$} [tr] at -11 -1 %
\arrow <2\einheit> [0.25,1] from -10 8.8 to -10 8.7 \plot -10 0  -10 20 /\put {$\onel$} [r] at -12.5 10 %
\put {\scriptsize $\bullet$} at -10 20 \put {$\zerotl$} [br] at -11 21 %
\put {\scriptsize $\bullet$} at 10 0 \put {$\zerobr$} [tl] at 11 -1 %
\arrow <2\einheit> [0.25,1] from 10 8.8 to 10 8.7 \plot 10 0  10 20 /\put {$\oner$} [l] at 12.5 10 %
\put {\scriptsize $\bullet$} at 10 20 \put {$\zerotr$} [bl] at 11 21 %
\arrow <2\einheit> [0.25,1] from 1.2 20 to 1.3 20 \plot -10 20  10 20 /\put {$\onet$} [b] at 0 22.5 %
\arrow <2\einheit> [0.25,1] from 1.2 0 to 1.3 0 \plot -10 0  10 0 /\put {$\oneb$} [t] at 0 -2.5 %
\put {$\twoc$} at 0 10 %
} %
}%

\def\leq{\leqslant}

\author{ R. Brown\footnote{Mathematics Department, Dean St.,
Bangor, Gwynedd LL57 1UT, UK. email: r.brown@bangor.ac.uk. Brown
has been supported for a part of this research by a Leverhulme
Emeritus Fellowship (2002-2004). }, K.H.
Kamps\footnote{Fachbereich Mathematik, FernUniversit\"at in Hagen,
D-58084 Hagen, Germany. email: heiner.kamps@fernuni-hagen.de}
  and T.Porter\footnote{Mathematics Department, Dean St.,
Bangor, Gwynedd LL57 1UT, UK. email: t.porter@bangor.ac.uk}}

\title{ A homotopy double groupoid of a Hausdorff
space II: \\ a van Kampen theorem}

{\theorembodyfont{\rmfamily}\newtheorem{example}{Example}[section]}
{\theorembodyfont{\rmfamily}\newtheorem{definition}[example]{Definition}}
\newtheorem{proposition}[example]{Proposition}
\newtheorem{theorem}[example]{Theorem}
{\theorembodyfont{\rmfamily}\newtheorem{rem}[example]{Remark}}

\newtheorem{lemma}[example]{Lemma}
\newcommand{\labto}[1]{\stackrel{#1}{\longrightarrow}}

\newcommand{\kreis}{\raisebox{0.5mm}[0mm][0mm]{\scriptsize$\circ$}}

\def\<{\langle}
\def\>{\rangle}

\def\bD{\mathbb{D}}
\def\bT{\mathbb{T}}

\newcommand{\tsq}{\mbox{\rule{0.04em}{1.55ex}\hspace{-0.00em}\rule{0.7em}{0.1ex}\hspace{-0.7em}\rule[1.5ex]{0.7em}{0.1ex}\hspace{-0.03em}\rule{0.04em}{1.55ex}}}

\newcommand{\br}{\mbox{\rule{0.7em}{0.2ex}\hspace{-0.04em}\rule{0.08em}{1.7ex}}}

\newcommand{\tl}{\mbox{\rule{0.08em}{1.7ex}\rule[1.54ex]{0.7em}{0.2ex}}}

\newcommand{\hh}{\mbox{\rule{0.7em}{0.2ex}\hspace{-0.7em}\rule[1.5ex]{0.70em}{0.2ex}}}

\newcommand{\vv}{\mbox{\rule{0.08em}{1.7ex}\hspace{0.6em}\rule{0.08em}{1.7ex}}}

\newcommand{\sq}{\mbox{\rule{0.08em}{1.7ex}\hspace{-0.00em}\rule{0.7em}{0.2ex}\hspace{-0.7em}\rule[1.54ex]{0.7em}{0.2ex}\hspace{-0.03em}\rule{0.08em}{1.7ex}}}

\def\brho{\boldsymbol{\rho}}

\def\hdgb{\boldsymbol{\rho}^\square}

\def\cU{\mathcal{U}}
\def\sD{\mathsf{D}}

\begin{document}
\maketitle \begin{center} {UWB Math Preprint 04.01} \end{center}
\begin{abstract}
This paper is the second in a series exploring the properties of a
functor which assigns a homotopy double groupoid with connections
to a Hausdorff space.

We show that this functor  satisfies a version of the van Kampen
theorem, and so is a suitable tool for nonabelian, 2-dimensional,
local-to-global problems. The methods are analogous to those
developed  by Brown and Higgins for similar theorems for other
higher homotopy groupoids.

An integral part of the proof is a  detailed discussion of
commutative cubes in a double category with connections, and a
proof of the key result that any composition of commutative cubes
is commutative. These results have recently been generalised to
all dimensions by Philip Higgins. \footnote{KEY WORDS: double
groupoid, double category, thin structure, connections,
commutative cube, van Kampen theorem
\\MATH SUBJECT CLASSIFICATION: 18D05,20L05,55Q05, 55Q35 }
\end{abstract}

\section*{Introduction}
A classical and key example of a nonabelian local-to-global
theorem  in dimension $1$ is the  van Kampen theorem for the
fundamental group of a space with base point:  if a space is the
union of two connected open sets with connected intersection, the
theorem determines  the fundamental group of the whole space, and
so the `global' information, in terms of the local information  on
the fundamental groups of the parts and the morphisms induced by
inclusions. Such a theorem thus relates a particular pushout of
spaces with base point to a pushout, or free product with
amalgamation, of their fundamental groups.

van Kampen's 1935 paper, \cite{vanKampen}, gave, in fact,  a
formula for the case of {\it non connected intersection}, as
required for the algebraic geometry applications he had in mind.
His formula follows from the version of the theorem  for the
fundamental groupoid of a space \cite{Brown1} (for the deduction
of the formula, see \S 8.4.9 of \cite{RBbook}). An outstanding
feature of this version is that it involves a universal property
which thus gives, in the non connected case, a {\it complete
determination of this nonabelian invariant} which includes the
fundamental groups based at various points. The basic reason for
this success seems to be that groupoids have structure in
dimensions 0 and 1, and so allow a better algebraic modelling of
the geometry of the intersections of the components of the two
open sets.

General theorems of this type, as in  for example \cite{colimits},
involve covers more general than by two open sets, and so require
the notion of {\it coequaliser} rather than pushout. Recall that
the following diagram is a basic example of a coequaliser in
topology:
\begin{equation} \label{coeqspace}
\bigsqcup _{(U,V)  \in \,           \mathcal{U}\, ^2}
 U\cap V \overset{a}{\underset{b}{\rightrightarrows}}
 \bigsqcup _{U \in\, \cU}U \labto{c}
 X  . \end{equation}Here $\mathcal U$ is an open cover
of the space $X$;  $\bigsqcup$ denotes coproduct in the category
of spaces,  which is given by disjoint union; the functions
$a,b,c$ are determined by the inclusions
$$a_{UV}:U\cap V\rightarrow U,\quad b_{UV}:U\cap V\rightarrow V,\quad c_U:U\rightarrow X$$ for
each $(U,V)  \in \mathcal{U} \,^2,U \in \cU$.  The coequaliser
condition  in this case says simply that any continuous function
$f:X\to Y$ is entirely determined by continuous functions $f_U: U
\to Y$ which agree on the intersections $U \cap V$. Thus a van
Kampen/coequaliser theorem states that the fundamental groupoid
functor applied to diagram \eqref{coeqspace} gives a coequaliser
diagram  of groupoids.

The paper \cite{HKK}, published in 2000,  showed the existence of
an {\it  absolute} homotopy 2-groupoid,
$\boldsymbol{\rho}^\bigcirc(X)$ of a space. By virtue of the
equivalence between 2-groupoids and  double groupoids with
connection proved in \cite{S,BM}, this also gave a homotopy double
groupoid with connection, $\hdgb(X)$, of a Hausdorff space $X$. An
explicit description of  $\hdgb(X)$ in terms of certain homotopy
classes of paths and squares, i.e. maps of $I$ and $I^2$, was
given in the first paper of this series, \cite{BHKP}. The interest
of these functors is that their 1-dimensional part
$\boldsymbol{\rho}_1$, which is the same for both functors,  has
the fundamental groupoid of $X$ as a quotient, and that they
contain the second homotopy groups $\pi_2(X,x)$ for all base
points $x$ in $X$.

We adapt the methods of \cite{BH} to prove a van
Kampen/coequaliser  theorem, i.e. a  2-dimensional local-to-global
property, for this functor $\hdgb$.

\noindent {\bf Theorem 4.1 [van Kampen/coequaliser theorem for
$\hdgb$]} {\it If $\cU$ is an open cover of the Hausdorff space
$X$, then the following diagram of morphisms induced by inclusions
\begin{equation}
\bigsqcup _{(U,V)  \in \,           \mathcal{U}\, ^2}
 \hdgb(U\cap V) \overset{a}{\underset{b}{\rightrightarrows}}
 \bigsqcup _{U \in\, \cU}\hdgb(U) \labto{c}
 \hdgb(X )  \label{coeq} \end{equation}
 is a coequaliser diagram in the category of double groupoids with
 connections.}

The morphisms      $a,b,c$ are determined by the inclusions
$a,b,c$ as above. Note also that  $\bigsqcup$, the coproduct in
the category of double groupoids with connections, is essentially
given by disjoint union.

The validity of such a theorem is an argument for the utility of
this particular algebraic structure. Note that the theorem yields
a new result in dimension 1, since the functor
$\boldsymbol{\rho}_1$ differs from the fundamental groupoid
$\pi_1$. In view of the success of the 1-dimensional fundamental
group(oid) in a variety of problems, as shown for example by a web
search on ``van Kampen theorem'', we hope that this theorem will
lead to wide investigations of this type of 2-dimensional
homotopical algebra.
\medskip

Here are some comments on the proof of the main theorem.

One key aspect of some proofs of a van Kampen type theorem is the
realisation of {\it algebraic inverse to subdivision}. In
dimension 1, this is more elegantly  realised by groupoids rather
than by groups. In dimension 2, this is more easily realised by
squares, using $\hdgb$, rather than by globes, using
$\boldsymbol{\rho}^\bigcirc$.

A second key aspect of a proof of the 1-dimensional van Kampen
theorem is the notion of {\it commutative square} in a groupoid,
and the fact that horizontal and vertical, and indeed any,
compositions of such commutative squares are commutative. In this
way, commutative squares in a groupoid form a double groupoid.
Likewise, a key aspect in the proof for higher dimensions is the
{\it definition of commutative cube} and the proof that any
composition of commutative cubes is commutative. That is, the
commutative cubes in a  double category should form a triple
category.

Our definition of commutative cubes requires the notion of
connections on a cubical set, as first introduced for double
groupoids in 1976 in \cite{BS:double}, and for double categories
in \cite{S}. We explain these ideas in the early sections of this
paper, and show their application to commutative cubes, here
called 3-shells. The proofs for compositions are deferred to a
penultimate section. These results have recently been generalised
to all dimensions by Philip Higgins in \cite{Higgins}, following
on from work in \cite{ABS,BM} and an earlier version of this
paper. We also note the earlier work of \cite{SW} on what they
call `homotopy commutative cubes'; this work  imposes more
conditions on the underlying double category than here, but also
uses connections.

The proof in section 4 of the main coequaliser theorem follows the
model for the 2-dimensional version for pairs of spaces given  in
\cite{BH}, with some extra twists. Given the main machinery, the
proof is quite short. It does not need the deformations based on
connectivity assumptions used in \cite{BH}. The proof in \cite{BH}
extends to a coequaliser theorem for the homotopy
$\omega$-groupoid of a filtered space, as in \cite{colimits}. This
suggests that the results of this paper should extend to higher
dimensions, but no construction is yet available of suitable {\it
absolute} higher homotopy groupoids.

In the final section, we explain intended future directions for
this work.

We would like to thank for helpful comments on an earlier draft:
Philip Higgins, who has  developed higher dimensional work in
\cite{Higgins}; and also Uli Fahrenberg, who has used some of
these ideas in directed homotopy theory, \cite{UF}.

\section{Double categories, connections and thin structures: definition and
notational conventions}

By a {\it double category} $ {\bf \sf K} $, we will mean what has
been called an {\it edge symmetric double category} in \cite{BM}.

We briefly recall some of the basic facts about double categories
in the above sense. In the first place, a double category, $ {\bf
\sf K} $, consists of a triple of category structures

\[
\begin{array}{c}
(K_{2} ,K_{1} ,\partial^{-}_{1} , \partial^{+}_{1} , +_{1} ,\varepsilon_{1} ) , \quad
(K_{2} , K_{1} ,\partial^{-}_{2} , \partial^{+}_{2} ,+_{2} ,\varepsilon_{2} )\\[2mm]
(K_{1} , K_{0} ,\partial^{-} , \partial^{+} , + , \varepsilon )
\end{array}
\]
 as partly shown in the diagram
\def\zoomfactor{1.2}\einheit=\zoomfactor mm%
\pictexdiagram / 2/ {%
\def\fontgroesse{10.95}\def\zeilenabstand{13.6}\fontsize{\fontgroesse}{\zeilenabstand}\selectfont%
\put {$K_1$} at 0 0 \put {$K_2$} at 0 20 \put {$K_0$} at 24 0 \put {$.$} at 28 -0.75 \put {$K_1$} at 24 20 %
\arrow <1.25\einheit> [0.25,1] from -1 17 to -1 3 \arrow <1.25\einheit> [0.25,1] from 1 17 to 1 3 %
\put {$\partial^-_1$} [r] at -2 10 \put {$\partial^+_1$} [l] at 2 10 %
\arrow <1.25\einheit> [0.25,1] from 23 17 to 23 3 \arrow <1.25\einheit> [0.25,1] from 25 17 to 25 3 %
\put {$\partial^-$} [r] at 22 10 \put {$\partial^+$} [l] at 26 10 %
\arrow <1.25\einheit> [0.25,1] from 4 -0.5 to 21 -0.5 \arrow <1.25\einheit> [0.25,1] from 4 1.5 to 21 1.5 %
\put {$\partial^+$} [t] at 12 -1.5 \put {$\partial^-$} [b] at 12 2.5 %
\arrow <1.25\einheit> [0.25,1] from 4 19.5 to 21 19.5 \arrow <1.25\einheit> [0.25,1] from 4 21.5 to 21 21.5 %
\put {$\partial^+_2$} [t] at 12 18.5 \put {$\partial^-_2$} [b] at 12 22.5 %
} %
The elements of $K_0,K_1,K_2$ will be called respectively {\it
points or objects, edges, squares}. The maps $ {\partial^\pm,} \;
\partial^{\pm}_{i} , \ i = 1,2 $, will be called {\it face maps},
the maps $ \varepsilon_{i} : K_{1} \longrightarrow K_{2} , \ i =
1,2 $, resp. $ \varepsilon : K_{0} \longrightarrow K_{1} $ will be
called {\it degeneracies}.

 The compositions, $ +_{1} $, resp. $ +_{2} $, are
referred to as {\it vertical} resp. {\it horizontal composition}
of squares. The axioms for a double category include the usual
relations of a 2-cubical set and the {\it interchange law}.
 {We use matrix notation for compositions as
$$\begin{bmatrix}a\\c \end{bmatrix} = a +_1 c, \quad
 \begin{bmatrix}a & b \end{bmatrix} = a +_2 b, $$ and the interchange law allows
one to use matrix notation
$$\begin{bmatrix} a&b \\c&d \end{bmatrix} $$
 for double composites of squares, as in \cite{BH2}. We also allow as in \cite{BH2}
 the multiple composition $[a_{ij}]$ of an array $(a_{ij})$ whenever for
 all appropriate $i,j$ we have
 $\partial ^+_1a_{ij}= \partial ^-_1a_{i+1,j},\;\partial
 ^+_2a_{ij}= \partial ^-_2a_{i,j+1}$.}

 The identities with respect to $ +_{1} $ ({\it vertical
identities}) are given by $ \varepsilon_{1} $ and will be
\def\zoomfactor{0.8}\einheit=\zoomfactor mm
denoted by \verticalidentity{0.6}\,. Similarly, we have {\it
horizontal identities} denoted by \horizontalidentity{0.6}\,.
Elements of the form $ \varepsilon_{1} \varepsilon (a) =
\varepsilon_{2} \varepsilon (a) $ for $ a \in K_{0} $ are called
{\it double degeneracies} and will be denoted by $ \odot $.

 A \textit{morphism of double categories} $f :
\mathsf{K}\to \mathsf{L}$ consists of a triple of maps $f_i : K_i
\to L_i$, $(i = 0,1,2)$, respecting the cubical structure,
compositions and identities.

 A {\it connection pair} on a double category $ { \mathsf K} $
 is given by a pair of maps
\[
\Gamma^{-} , \Gamma^{+} : K_{1} \longrightarrow K_{2}
\]
 whose edges are given by the following diagrams for $ a
\in K_{1} $:
\def\zoomfactor{0.8}\einheit=\zoomfactor mm%
\pictexdiagram / / {%
\def\fontgroesse{10.95}\def\zeilenabstand{13.6}\fontsize{\fontgroesse}{\zeilenabstand}\selectfont%
\put {$\Gamma^-(a)$} at -50 10 \put {$=$} at -40 10 %
\arrow <2\einheit> [0.25,1] from -30 8.8 to -30 8.7 \plot -30 0  -30 20 /\put {$a$} [r] at -32.5 10 %
\arrow <2\einheit> [0.25,1] from -10 8.8 to -10 8.7 \plot -10 0  -10 20 /\put {$1$} [l] at -7.5 10 %
\setquadratic\multiput {\vdeg} at -10 1 *9 0 2 /\setlinear%
\arrow <2\einheit> [0.25,1] from -18.8 20 to -18.7 20 \plot -30 20  -10 20 /\put {$a$} [b] at -20 22.5 %
\arrow <2\einheit> [0.25,1] from -18.8 0 to -18.7 0 \plot -30 0  -10 0 /\put {$1$} [t] at -20 -2.5 %
\setquadratic\multiput {\hdeg} at -30 0 *9 2 0 /\setlinear%
\put {$=$} at 0 10 %
\arrow <2\einheit> [0.25,1] from 10 8.8 to 10 8.7 \plot 10 0  10 20 /\put {$$} [r] at 7.5 10 %
\arrow <2\einheit> [0.25,1] from 30 8.8 to 30 8.7 \plot 30 0  30 20 /\put {$$} [l] at 32.5 10 %
\setquadratic\multiput {\vdeg} at 30 1 *9 0 2 /\setlinear%
\arrow <2\einheit> [0.25,1] from 21.2 20 to 21.3 20 \plot 10 20  30 20 /\put {$$} [b] at 20 22.5 %
\arrow <2\einheit> [0.25,1] from 21.2 0 to 21.3 0 \plot 10 0  30 0 /\put {$$} [t] at 20 -2.5 %
\setquadratic\multiput {\hdeg} at 10 0 *9 2 0 /\setlinear%
\put {$=$} at 40 10 %
\linethickness1pt\linethickness=\zoomfactor\linethickness%
\put {$\gammaminus{}$} at 20 10 \put {$\gammaminus{}$} at 50 10 %
\arrow <1.75\einheit> [0.25,1] from 70 20 to 70 11.25 \arrow <1.75\einheit> [0.25,1] from 70 20 to 78.75 20 %
\put {\footnotesize $1$} [t] at 70 9.25 \put {\footnotesize $2$} [l] at 80.25 20 %
} %

\def\zoomfactor{0.8}\einheit=\zoomfactor mm%
\pictexdiagram / -8.75/ {%
\def\fontgroesse{10.95}\def\zeilenabstand{13.6}\fontsize{\fontgroesse}{\zeilenabstand}\selectfont%
\put {$\Gamma^+(a)$} at -50 10 \put {$=$} at -40 10 %
\arrow <2\einheit> [0.25,1] from -30 8.8 to -30 8.7 \plot -30 0  -30 20 /\put {$1$} [r] at -32.5 10 %
\setquadratic\multiput {\vdeg} at -30 1 *9 0 2 /\setlinear%
\arrow <2\einheit> [0.25,1] from -10 8.8 to -10 8.7 \plot -10 0  -10 20 /\put {$a$} [l] at -7.5 10 %
\arrow <2\einheit> [0.25,1] from -18.8 20 to -18.7 20 \plot -30 20  -10 20 /\put {$1$} [b] at -20 22.5 %
\setquadratic\multiput {\hdeg} at -30 20 *9 2 0 /\setlinear%
\arrow <2\einheit> [0.25,1] from -18.8 0 to -18.7 0 \plot -30 0  -10 0 /\put {$a$} [t] at -20 -2.5 %
\put {$=$} at 0 10 %
\arrow <2\einheit> [0.25,1] from 10 8.8 to 10 8.7 \plot 10 0  10 20 /\put {$$} [r] at 7.5 10 %
\setquadratic\multiput {\vdeg} at 10 1 *9 0 2 /\setlinear%
\arrow <2\einheit> [0.25,1] from 30 8.8 to 30 8.7 \plot 30 0  30 20 /\put {$$} [l] at 32.5 10 %
\arrow <2\einheit> [0.25,1] from 21.2 20 to 21.3 20 \plot 10 20  30 20 /\put {$$} [b] at 20 22.5 %
\setquadratic\multiput {\hdeg} at 10 20 *9 2 0 /\setlinear%
\arrow <2\einheit> [0.25,1] from 21.2 0 to 21.3 0 \plot 10 0  30 0 /\put {$$} [t] at 20 -2.5 %
\put {$=$} at 40 10 %
\linethickness1pt\linethickness=\zoomfactor\linethickness%
\put {$\gammaplus{}$} at 20 10 \put {$\gammaplus{}$} at 50 10 %
\put {.} at 60 0 %
} %

 The axioms for the connection pair are listed in
\cite{BM}, Section 4. In particular, the {\it transport laws}
\begin{equation} \label{transport} \begin{bmatrix}\; \tl & \hh \;\\ \vv & \tl \end{bmatrix} \; =
\; \tl\;, \qquad \begin{bmatrix}\; \br & \vv\; \\ \hh & \br
\end{bmatrix} \; = \; \br\end{equation}
describe the connections of the composition of two elements, while
the laws
\begin{equation}\label{cancel} \begin{bmatrix}\; \tl \; \\ \br \end{bmatrix}\;=\; \hh\;, \qquad \begin{bmatrix}\; \tl &  \br \;\end{bmatrix}\;=\;
\vv,\end{equation} allow cancellation of connections.

 A {\it morphism of double categories with
connections} is a morphism of double categories respecting
connections.

 Brown and Mosa in \cite{BM} have shown that a pair of
connections on a double category $\mathsf{K}$ is equivalent to a
{\it  thin structure} on $\mathsf{K}$, whose definition we recall.

  We recall the standard language for `commutative squares'; we
  will later move to `commutative cubes'.

 First, if $\alpha$ is a square in $\mathsf{K}$, then the {\it
boundary $(2$-shell$)$ } of $\alpha $ is the quadruple
$$(\partial^-_2 \alpha,
\partial^+_1 \alpha,\partial^-_1 \alpha,\partial^+_2 \alpha).$$
We also say $\alpha$ is a {\it filler} of its boundary. This
boundary {\it commutes} if $\partial^-_2 \alpha+\partial^+_1
\alpha=\partial^-_1 \alpha+\partial^+_2 \alpha$. More generally, a
2-{\it shell} is defined to be a quadruple $(a,b,c,d)$ of edges
such that $\partial^- a=\partial^-c,
\partial^+ a = \partial^- b, \partial^+c=\partial^- d, \partial^+b = \partial^+ d$,
and this 2-shell {\it commutes}  if $a+b=c+d$.

If $C$ is a category, then by $\tsq C$ we denote the double
category of commuting squares (2-shells) in $C$ with the obvious
double category structure.  Then a {\it thin structure} $\Theta$
on a double category $\mathsf{K}$ is a morphism of double
categories
$$\Theta : \tsq K_1 \to \mathsf{K},$$
which is the identity on  {$K_1$.}  The elements of $K_2$ lying in
$\Theta(\tsq K_1)$ are called {\it thin  squares}. The definition
implies immediately:

 \vspace{0.5ex}
 (T0) \emph{A thin  square has commuting boundary.}

 \vspace{0.5ex}

(T1) \emph{Any commutative  {$2$-shell} has a unique thin filler.}

\vspace{0.5ex} (T2) \emph{{Any identity square, }  and any
composition of thin squares,  is thin.}

\medskip

In the thin structure induced on a double category by a pair of
connections, any thin square is a certain composite of identities,
both horizontal and vertical, and connections. For an explicit
formula, we refer to \cite{BM}, Theorem 4.3,  and we will refer to
these thin squares as being \emph{algebraically thin}. Of course,
a morphism of double categories with connections will preserve
(algebraic) thinness.

\section{Commutative cubes}

We need to extend the domain of discourse from double categories
to triple categories in order to explain the notion of {\it
commutative cube  {($3$-shell)} in a double category with
connections}. This notion has been defined  {for $n$-shells} in
the case of cubical $\omega$-groupoids with connections in
\cite[Section 5]{BH2}, and for the more difficult category case in
\cite[Section 9]{ABS}. For the convenience of the reader we set up
the theory in our low dimensional case.

Let $\bD, \; \bT$ be respectively the categories of double
 {categories}  and of triple  {categories} with connections,
in the sense of \cite{ABS}. There is a natural and obvious
truncation functor $tr: \bT \to \bD$, which forgets the
3-dimensional structure. This functor has a right adjoint $cosk:
\bD \to \bT$ which may be constructed  in terms of {\it cubes} or
3-{\it shells}, which we now define.

\begin{definition} Let $ { \mathsf K} $ be a double
category. A {\it cube}  { (\emph{$3$-shell})} in $ { \mathsf K} $,
\[
\alpha = (\alpha^{-}_{1} , \alpha^{+}_{1} , \alpha^{-}_{2} , \alpha^{+}_{2} , \alpha^{-}_{3} ,
\alpha^{+}_{3} )
\]
 consists of squares $ \alpha^{\pm}_{i} \in K_{2} \quad (i = 1,2,3) $ such that
\[
\partial^{\sigma}_{i} ( \alpha^{\tau}_{j} ) = \partial^{\tau}_{j-1} (\alpha^{\sigma}_{i} )
\]
 for $ \sigma , \tau = \pm  $ and $ 1 \leq i < j \leq 3 $. \hfill
 $\blacksquare$
\end{definition}

 A cube may be illustrated by the following picture:

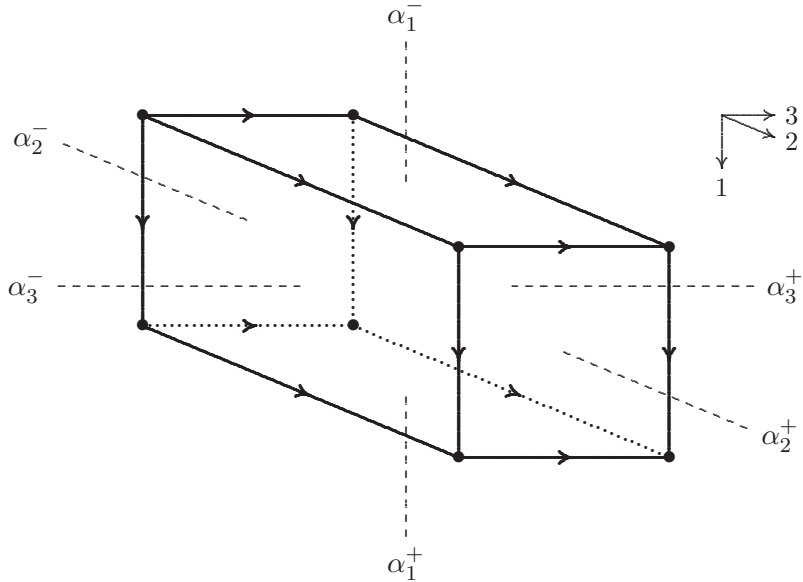
\backsquarelabels &&& %
\middlearrowlabels &&& %
\frontsquarelabels &&& %
\facelabels \alpha^-_2&\alpha^+_2&\alpha^-_3&\alpha^+_3&\alpha^+_1&\alpha^-_1 %
\def\zoomfactor{0.7}\def\fontgroesse{10.95}\def\zeilenabstand{13.6}\def\showindexing{1}%
\letteredcube &&& %

\begin{definition} \label{compcubes}
 We define  three partial compositions, $ +_{1} , +_{2} ,
+_{3} $, of cubes.  Let $ \alpha , \beta $ be cubes in $ {\bf \sf
K} $.
\begin{enumerate}[(i)]
\item If  $ \alpha^{+}_{1} = \beta^{-}_{1} $, then we define
           \[
           \alpha +_{1} \beta = (\alpha^{-}_{1} , \beta^{+}_{1} , \alpha^{-}_{2} +_{1} \beta^{-}_{2} ,
           \alpha^{+}_{2} +_{1} \beta^{+}_{2} , \alpha^{-}_{3} +_{1} \beta^{-}_{3} , \alpha^{+}_{3} +_{1}
           \beta^{+}_{3} ).
           \]
\item If  $ \alpha^{+}_{2} = \beta^{-}_{2} $, then we define
           \[
           \alpha +_{2} \beta = (\alpha^{-}_{1} +_{1} \beta^{-}_{1} , \alpha^{+}_{1} +_{1} \beta^{+}_{1} ,
           \alpha^{-}_{2} , \beta^{+}_{2} , \alpha^{-}_{3} +_{2} \beta^{-}_{3} , \alpha^{+}_{3} +_{2}
           \beta^{+}_{3} ).
           \]
\item If  $ \alpha^{+}_{3} = \beta^{-}_{3} $, then we define
           \[
           \alpha +_{3} \beta = ( \alpha^{-}_{1} +_{2} \beta^{-}_{1} , \alpha^{+}_{1} +_{2}
           \beta^{+}_{1} , \alpha^{-}_{2} +_{2} \beta^{-}_{2} , \alpha^{+}_{2} +_{2}
           \beta^{+}_{2} , \alpha^{-}_{3} , \beta^{+}_{3} ).
           \]
\end{enumerate}
\end{definition}

This is a special case of general definitions in dimension $n$
given in \cite{BH2}, see also \cite{Higgins}.

As explained earlier, we need the notion of a 3-shell $\alpha$
being `commutative'. Intuitively, this says that the composition
of the odd faces $\alpha^-_1,\alpha^-_3,\alpha^+_2$ of $\alpha$ is
equal to the composition of the even faces
$\alpha^+_1,\alpha^+_3,\alpha^-_2$ (think of -,+ as 0,1). The
problem of how to make sense of such compositions of three squares
in a double groupoid, and obtain the right boundaries, is solved
using the connections.

{\begin{definition} \label{hal} Suppose given, { in a double
category with connections $ {\mathsf K} $, a cube (3-shell)}
$$
\alpha = (\alpha^{-}_{1} , \alpha^{+}_{1} , \alpha^{-}_{2} ,
\alpha^{+}_{2} , \alpha^{-}_{3} , \alpha^{+}_{3} ).
$$
 We define the {\it composition of the odd faces} of $\alpha$  to
be
\begin{align} \boldsymbol{\partial}^{\mathrm{odd}} \alpha &=
\begin{bmatrix}\;  \tl & \alpha^-_1 & \br \\ \alpha^-_3 & \alpha^+_2
& \hh
\end{bmatrix} \\
\intertext{and the {\it composition of the even faces} of $\alpha$
to be}\boldsymbol{\partial}^{\mathrm{even}} \alpha &=
\begin{bmatrix}\hh & \alpha^-_2 & \alpha^+_3 \\
\tl & \alpha^+_1 & \br\;
\end{bmatrix}
\end{align}
This definition can be regarded as a cubical, categorical (rather
than groupoid) form of the Homotopy Addition Lemma (HAL) in
dimension 3. We refer to the discussion below
\end{definition} }
\begin{definition} \label{commcube}

We define $\alpha$ to be {\it commutative} if it satisfies the
Homotopy Commutativity Lemma (HCL), i.e. \begin{equation}
\boldsymbol{\partial}^{\mathrm{odd}} \alpha=
\boldsymbol{\partial}^{\mathrm{even}} \alpha. \tag{HCL}
\end{equation}
\end{definition}

The reader should draw a 3-shell, label all the edges with
letters, and see that this equation makes sense in that the
2-shells of each side of equation (HCL) coincide. Notice however
that these 2-shells do not have coincident partitions along the
edges: that is the edges of this 2-shell in direction 1  are
formed from different compositions of the type $0+a$ and $a +0$.
{This formula should also be compared with the notion of `homotopy
commutative cube' defined in a more restricted kind of double
category with connections in \cite{SW}.}

The formula in Definition \ref{commcube} is unsymmetrical, and
this seems in part a consequence of trying to express a
3-dimensional idea in a 2-dimensional formula. It also reflects
the difficulty of the concept. A formula in dimension 4 is given
in \cite{Gaucher}.

There are other  pictures and forms  of this which we now explore,
and some of which will be used later.

First we give a picture in another format which also shows how the
boundary of each side is made up:

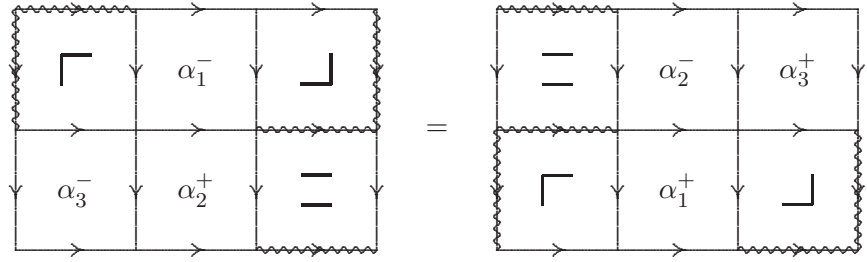
\onecubeleft &&&&& %
\twocubeleft \alpha^-_3&\alpha^+_2&\alpha^-_1 %
\onecuberight &&&&& %
\twocuberight \alpha^+_1&\alpha^-_2&\alpha^+_3 %
\def\zoomfactor{0.8}\def\fontgroesse{10.95}\def\zeilenabstand{13.6}\def\orientation{0}\def\showindexing{0}%
\letteredrectangles &&&&=

This is in fact equivalent to:

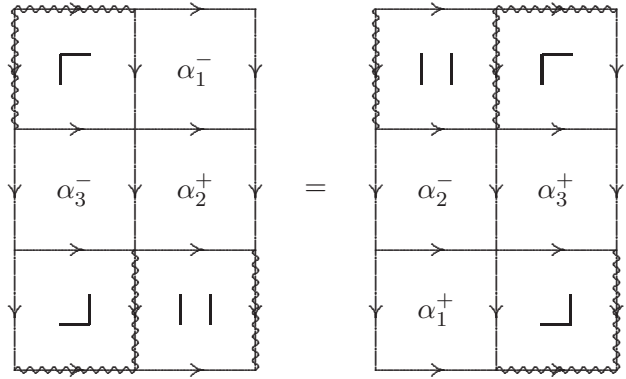
\onecubeleft &&&&& %
\twocubeleft \alpha^-_3&\alpha^+_2&\alpha^-_1 %
\onecuberight &&&&& %
\twocuberight \hspace*{28mm}\phantom{.}\alpha^+_1\hspace*{28mm}.&\alpha^-_2&\alpha^+_3 %
\def\zoomfactor{0.8}\def\fontgroesse{10.95}\def\zeilenabstand{13.6}\def\orientation{1}\def\showindexing{0}%
\letteredrectangles &&2.75&&= %

 The reader can check that in each case the two sides of the
equation have the same 2-shells (but not as subdivided).

The second equation can also be written  in matrix notation as: 
\begin{align*}
\left [\begin{array}{cc}\gammaplus{0.6}\rule{0mm}{7\einheit}&\alpha^-_1\\[2.5mm]
\alpha^-_3&\alpha^+_2\\[2.5mm]
\gammaminus{0.6}&\verticalidentity{0.6}\end{array}\right ]& =
\left [\begin{array}{cc}\verticalidentity{0.6}&\gammaplus{0.6}\rule{0mm}{7\einheit}\\[2.5mm]
\alpha^-_2&\alpha^+_3\\[2.5mm]
\alpha^+_1&\gammaminus{0.6}\end{array}\right ].
\tag{HCL$'$}\end{align*}

 {The equivalence of (HCL)  to (HCL$'$)  can be
seen by adding $\begin{bmatrix} \odot & \verticalidentity{0.6}
&\gammaplus{0.6}\;
\end{bmatrix}$, (respectively $\begin{bmatrix}\gammaminus{0.6} & \verticalidentity{0.6}&
 \odot
\end{bmatrix}$) to the top (resp. bottom) of the two sides of (HCL), using equations
\eqref{cancel} for further cancellation of
connections, and then absorbing identities.}

 We note that a morphism of double categories with
connections preserves  commutativity of cubes.

 The following theorem will be crucial for the proof of the van Kampen theorem.

 \begin{theorem}\label{compcubes}
 Let $ {\bf \sf K} $ be a double category with connections. Then
any composition of commutative $3$-shells is commutative.
\end{theorem}
The proof is left to section \ref{comproof}.

\begin{proposition}
The functor $tr: \bT \to \bD$ has a left adjoint $sk: \bD \to \bT$
which assigns to a double groupoid with connections {\bf \sf D}
the triple groupoid $\sq\, ${\bf \sf D} which agrees with {\bf \sf
D} in dimensions $\leq 2$ and in dimension $3$  consists of the
commutative $3$-shells in {\bf \sf D}.
\end{proposition}

\begin{rem}  (i) Our  homotopy commutativity lemma has been called
also the `homotopy addition lemma', and has been formulated in
different equivalent forms in certain $ 2 $-dimensional groupoid
type settings using inverses and reflections (cf. \cite{BH},
Proposition 5, \cite{SW}, Definition 3.8, \cite{BHKP}, Proposition
5.5).  Also \cite{BH2}, \S 5, 7, sets up, in the situation of
cubical multiple groupoids with connections, truncation and
skeleton functors tr$^n$ and sk$^n$ and an $n$-dimensional
homotopy addition lemma as a formula for the boundary of an
$n$-cube. The definition and use of thin elements in dimension $n$
is a key part of the proof of the generalised van Kampen theorem
in \cite{colimits}.

The homotopy addition lemma presented here applies to a general
2-dimensional \emph{category} type setting. It is also a special
case of results on cubical $\omega$-categories with connections in
\cite[Section 9]{ABS}. (The fact that any composition of
commutative shells is commutative is not stated in the last
reference, but does follow easily from the results stated there.)

(ii) Homotopy commutative cubes have been defined and investigated
by Spencer and Wong in the more special setting of what they call
{\it special double categories with connection} (\cite{SW},
Definition 1.4). By this they mean a double category (with
connections), $ {\bf \sf K} $, such that the horizontal subdouble
category of squares in $ {\bf \sf K} $ of the form 
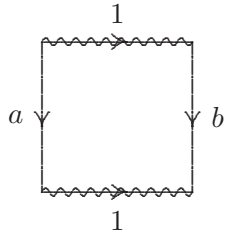
\def\zoomfactor{1}\einheit=\zoomfactor mm%
\pictexdiagram / / {%
\def\fontgroesse{10.95}\def\zeilenabstand{13.6}\fontsize{\fontgroesse}{\zeilenabstand}\selectfont%
\arrow <2\einheit> [0.25,1] from -10 8.8 to -10 8.7 \plot -10 0  -10 20 /\put {$a$} [r] at -12.5 10 %
\arrow <2\einheit> [0.25,1] from 10 8.8 to 10 8.7 \plot 10 0  10 20 /\put {$b$} [l] at 12.5 10 %
\arrow <2\einheit> [0.25,1] from 1.2 20 to 1.3 20 \plot -10 20  10 20 /\put {$1$} [b] at 0 22.5 %
\setquadratic\multiput {\hdeg} at -10 20 *9 2 0 /\setlinear%
\arrow <2\einheit> [0.25,1] from 1.2 0 to 1.3 0 \plot -10 0  10 0 /\put {$1$} [t] at 0 -2.5 %
\setquadratic\multiput {\hdeg} at -10 0 *9 2 0 /\setlinear%
} %
 is a groupoid under $ +_{2} $. Spencer and Wong  {show} that
 this is a suitable setting to deal with (homotopy)
pullbacks/pushouts and gluing/cogluing theorems in homotopy
theory.

It is not difficult to show that the definition of a homotopy
commutative cube as given in \cite{SW}, Definition 3.8, is
equivalent to our definition of commutative cube.
\end{rem}

 As we mentioned before the fact that
commutativity of cubes is preserved under composition will be
crucial for our purposes. Spencer and Wong have a similar result
in their setting (\cite{SW}, Proposition 3.11).

The basic notions of relative homotopy carry over from topology to
algebra. Let $\mathsf{K}$ be a double category with connections. A
square $u$ in $\mathsf{K}$ is called a {\it relative homotopy} if
the two opposite edges  {$\partial^-_2 u,
\partial^+_2 u$} are identities. A {\it relative homotopy between
squares} $u$, $u^\prime$ in $\mathsf{K}$ is a commutative cube
$\alpha$ such that
$$\alpha_1^- = u, \quad \alpha_1^+= u^\prime$$
and the remaining faces are relative homotopies.  If, in addition,
the remaining faces are thin, then $u$ and $u^\prime$ are called
{\it thinly equivalent}

 From (T0) and (T1), we conclude:

\medskip
(T3) \emph{A thin square which is a relative homotopy is an identity.}
\medskip

 From Definition 2.2, making use of \cite{BM}, 3.1(i), we read off
 the following  { lemma:}

 \begin{lemma} \label{lem:coincide} Squares which are thinly equivalent, coincide.
 \end{lemma}
 {This key lemma is used later to show that a construction is
independent of the choices made. }

\noindent
\section{The homotopy double groupoid,  $\boldsymbol{\rho}^{\square} (X) $}
This section is adapted from  {\cite{BHKP}}, and the reader should
refer to that source for fuller details.
\subsection{The singular cubical set of a topological space}
We shall be concerned with the low dimensional part (up to
dimension 3) of the singular cubical set  $$R^{\square}
(X)=(R^{\square}_n(X),\partial^{-}_i, \partial^{+}_{i} ,
\varepsilon_i)$$ of a topological space $X.$ We recall the
definition (cf. \cite{BH2}).

For $n\geqslant 0$ let
$$R^{\square}_n(X)= \ \mbox{\bf Top} (I^n,X)$$
denote the set of {\it singular} $n$--{\it cubes} in $X,$ i.e. continuous maps
$I^n\longrightarrow X,$ where $I=[0,1]$ is the unit interval of real
numbers.

We shall identify $R^{\square}_0(X)$ with the set of points of $X.$
For $n=1,2,3$ a singular $n$--cube
will be called a {\it path}, resp. {\it square}, resp.
{\it cube}, in $X.$

The {\it face maps}  $$\partial^{-}_i ,\partial^{+}_{i}
:R^{\square}_n (X)\longrightarrow R^{\square}_{n-1} (X)\quad
(i=1,\ldots ,n)$$ are given by inserting $ 0 $ resp. $ 1 $ at the
$i^{th} $ coordinate whereas the {\it degeneracy maps}  $$
\varepsilon_{i} : R^{\square}_{n-1} (x) \longrightarrow
R^{\square}_{n} (X) \quad (i=1,\ldots ,n) $$ are given by omitting
the $ i^{th} $ coordinate. The face and degeneracy maps satisfy
the usual cubical relations (cf. \cite{BH2}, {\S} 1.1; \cite{KP} ,
{\S} 5.1). A path $ a \in R^{\square}_{1} (X) $ has {\it initial
point} $ a(0) $ and {\it endpoint} $ a(1) $. We will use the
notation $ a : a(0) \simeq a(1) $. If $ a,b $ are paths such that
$ a(1) = b(0) $, then we denote by $ a+b : a(0) \simeq b(1) $
their \emph{concatenation}, i.e.
$$ (a+ b)(s)= \left\{\begin{array}{ll}
a(2s),& 0\leqslant
s\leqslant\frac{1}{2}\\
       b(2s-1),& \frac{1}{2}\leqslant s\leqslant 1\end{array}\right.$$
If $x$ is a point of $X,$ then $\varepsilon_1 (x)\in
R^{\square}_1(X),$ denoted $e_x,$ is the \emph{constant path} at
$x,$ i.e.  $$e_x(s)=x \mbox{ for all } s\in I.$$

If $ a : x \simeq y $ is a path in $ X $, we denote by $ -a : y
\simeq x $ the {\it path reverse} to  $ a $, i.e.  $ (-a)(s) =
a(1-s) $ for $ s \in I $. In the set of squares $R^{\square}_2(X)$
we have two partial compositions $+_{1}$ (\emph{vertical
composition}) and $+_{2}$ (\emph{ horizontal composition}) given
by concatenation in the first resp. second variable.

 Similarly, in
the set of cubes $ R^{\square}_{3} (X) $ we have three partial
compositions $ +_{1} , +_{2} , +_{3} $.

The standard properties of vertical and horizontal composition of
squares are listed in \cite{BH2},  {\S}1. In particular we have
the following \emph{interchange law}. Let $u,u',\,w,w'\in
R^{\square}_2(X)$ be squares, then
$$ (u +_{2} w)+_{1} (u' +_{2} w')=(u+_{1}
u')+_{2} (w+_{1} w'),$$whenever both sides are defined. More
generally, we have an interchange law for rectangular
decomposition of squares. In more detail, for positive integers
$m,n$ let $\varphi_{m,n}:I^2\longrightarrow [0,m]\times [0,n]$ be
the homeomorphism $(s,t)\longmapsto (ms,nt).$ An $m\times n$ {\it
subdivision} of a square $u:I^2\longrightarrow X$ is a
factorization $u=u'\,\kreis\, \varphi_{m,n};$ its {\it parts} are
the squares $u_{ij}: {I^2}\longrightarrow X$ defined by
$$u_{ij}(s,t)=u'(s+i-1,t+j-1).$$
We then say that $u$ is the \emph{composite} of the  {array of
squares $(u_{ij})$}, and we use matrix notation $u=[u_{ij}].$ Note
that as in \S 1, $u+_{1} u',\enskip u +_{2} w$ and the two sides
of the interchange law can be written respectively as
$$\begin{bmatrix}u \\ u'\end{bmatrix}, \qquad [u\enskip w], \qquad
\begin{bmatrix}u & w\\
               u' & w'\end{bmatrix}\,.$$

Finally, \emph{ connections}
$$\Gamma^{-} ,\Gamma^{+} : R^{\square}_{1} (X) \longrightarrow R^{\square}_{2} (X) $$
are defined as follows.  {If}  $ a \in R^{\square}_{1} (X) $ is a
path, $ a:x \simeq y $, then let
\[
\Gamma^{-} (a)(s,t) = a(\mbox{max} (s,t)); \ \Gamma^{+} (a)(s,t) = a
(\mbox{min} (s,t)).
\]
The full structure of $ R^{\square} (X) $ as a {\it cubical
complex with connections and compositions} has been exhibited in
 {\cite{ABS,grandis-mauri}.}

\subsection{Thin squares}
In the setting of a geometrically defined double groupoid with
connection, as in \cite{BH2}, (resp \cite{BHKP}), there is an
appropriate notion of \emph{geometrically thin} square.  {It is
proved in \cite{BH2}, Theorem 5.2
  (resp. \cite{BHKP}, Proposition 4), that in the  cases given there}, geometrically and
algebraically thin squares coincide. In our context the explicit definition is as follows:

\begin{definition}\label{thin} (1) A square $ u:I^{2} \longrightarrow X $ in a topological
space $ X $ is \emph{thin} if there is a factorisation of $ u $
$$ u : I^{2} \stackrel{\Phi_{u}}{\longrightarrow}
J_{u} \stackrel{p_{u}}{\longrightarrow} X ,$$ where $ J_{u} $ is a
tree and $ \Phi_{u} $ is piecewise linear (PWL, see below) on the
boundary $ \partial{I}^{2} $ of $ I^{2} $.

Here, by a {\it tree}, we mean the underlying space $ |K| $ of a
finite $ 1 $-connected $ 1 $-dimensional simplicial complex $ K $.

A map $ \Phi : |K| \longrightarrow |L| $ where $ K $ and $ L $ are
(finite) simplicial complexes is PWL ({\it piecewise linear}) if
there exist subdivisions of $ K $ and $ L $ relative to which $
\Phi $ is simplicial.

(2) Let $ u$ be as above,  then the homotopy class of $ u $
relative to the boundary $ \partial{I}^{2} $ of $ I $ is called a
{\it double track}.  A double track is {\it thin} if it has a thin
representative. \hfill $\blacksquare$
\end{definition}
{We note that the class of thin squares  in a topological space
$X$ is closed under vertical and horizontal composition of squares
if $X$ is assumed to be Hausdorff (\cite[Proposition 3.5]{BHKP}).
This is a consequence of appropriate pushout properties of trees
(cf. \cite[Section 5]{BHKP}).}

\subsection{The homotopy double groupoid of a Hausdorff space}

The data for the homotopy double groupoid, $ \brho^{\square} (X)
$, {of a Hausdorff space $X$} will be denoted by
\[
\begin{array}{c}
(\brho^{\square}_{2} (X) , \brho^\square_1 (X) ,
\partial^{-}_{1} , \partial^{+}_{1} , +_{1} , \varepsilon_{1} ) ,
\ (\brho^{\square}_{2} (X), \brho^\square_1 (X) ,
\partial^{-}_{2} , \partial^{+}_{2} , +_{2} , \varepsilon_{2} )\\[3mm]
(\brho^\square_1 (X) , X , \partial^{-} ,
\partial^{+} , + , \varepsilon ).
\end{array}
\]
 Here $\brho^\square_1 (X)$ denotes the \emph{path groupoid} of $X$
of \cite{HKK}. We recall the definition. The objects of $
\brho^\square_1 (X) $ are the points of $ X $. The morphisms of $
\brho^\square_1 (X) $ are the equivalence classes of paths in $ X
$ with respect to the following relation $ \sim_{T} $.
\medskip

\begin{definition}
Let $ a,a' : x \simeq y $ be paths in $ X $. Then
$ a$ is \emph{ thinly equivalent} to $ a' $, denoted $ a \sim_{T} a' $, if
there is a thin relative homotopy between $ a $ and $ a' $.
\end{definition}

We note that $ \sim_{T} $ is an equivalence relation, see
\cite{BHKP}. We use $ \langle a \rangle : x \simeq y $ to denote
the $ \sim_{T} $ class of a path $ a: x \simeq y $ and call $
\langle a \rangle $  the {\it semitrack} of $ a $. The groupoid
structure of $ \brho^\square_1 (X) $ is induced by concatenation,
+, of paths. Here one makes use of the fact that if $ a: x \simeq
x', \ a' : x' \simeq x'', \ a'' : x'' \simeq x''' $ are paths then
there are canonical thin relative homotopies
\[
\begin{array}{r}
(a+a') + a'' \simeq a+ (a' +a'') : x \simeq x''' \ ({\it rescale}) \\
a+e_{x'} \simeq a:x \simeq x' ; \ e_{x} + a \simeq a: x \simeq x' \
({\it dilation}) \\
a+(-a) \simeq e_{x} : x \simeq x \ ({\it cancellation}).
\end{array}
\]
The source and target maps of  $\brho^\square_1 (X)$ are given by
$$\partial^{-}_{1} \langle a\rangle =x,\enskip \partial^{+}_{1}
\langle a\rangle =y,$$
if $\langle a\rangle :x\simeq y$ is a semitrack. Identities and inverses
are given by
$$\varepsilon (x)=\langle e_x\rangle  \quad \mathrm{ resp.} -\langle a\rangle
=\langle -a \rangle.$$
In order to construct $\brho^{\square}_2 (X)$, we define a
relation of cubically thin homotopy on the set $R^{\square}_2(X)$
of squares.

\begin{definition}  Let $u,u'$ be squares in $X$ with common vertices.
(1) A {\it cubically thin homotopy} $U:u\equiv^{\square}_T u'$
between $u$ and $u'$ is a cube $U\in R^{\square}_3(X)$ such that

(i) $U$ is a homotopy between $u$ and $u',$

\begin{center}
i.e. $\partial^{-}_1 (U)=u,\enskip \partial^{+}_1 (U)=u',$\end{center}
(ii) $U$ is rel. vertices of $I^2,$
\begin{center}
i.e. $\partial^{-}_2\partial^{-}_2 (U),\enskip\partial^{-}_2
\partial^{+}_2 (U),\enskip
\partial^{+}_2\partial^{-}_2 (U),\enskip\partial^{+}_2
\partial^{+}_2 (U)$ are
constant,\end{center} (iii) the faces $ \partial^{\alpha}_{i} (U)
$ are thin for {$ \alpha = \pm 1, \ i = 2,3 $}.

(2) The square $u$ is {\it cubically} $T$-{\it equivalent} to
$u',$ denoted $u\equiv^{\square}_T u'$ if there is a cubically
thin homotopy between $u$ and $u'.$
\end{definition}

\medskip

\begin{proposition}
The relation $\equiv^{\square}_T$  is
an equivalence relation on $R^{\square}_2(X).$
\end{proposition}
\textbf{Proof} The reader is referred to \cite{BHKP} for a proof.
\hfill$\blacksquare$

\medskip

If $u\in R^{\square}_2(X)$ we write $ \{ u \}^{\square}_{T} $, or
simply $\{u\}_T$, for the equivalence class of $u$ with respect to
$\equiv^{\square}_T.$ We denote the set of equivalence classes
$R^{\square}_2(X)\equiv^{\square}_T$ by $\brho^{\square}_2 (X).$
This inherits the operations and the geometrically defined
connections from $R^{\square}_2(X)$ and so becomes a double
groupoid with connections.  A proof of the final fine detail of
the structure is given in \cite{BHKP}.

\medskip

\begin{definition}
An element of $\brho^{\square}_2 (X)$  is \emph{thin} if it has a
thin representative (in the sense of Definition \ref{thin}).
\end{definition}

From the remark at the beginning of this subsection we infer:

\begin{lemma}
Let $f: \brho^\square(X) \to \mathsf D$ be a morphism of double
groupoids
 with connection.  {If $\alpha \in {\brho^\square_2}(X)$ is thin,
 then $f(\alpha)$ is thin.}\hfill$\blacksquare$
\end{lemma}
\subsection{The homotopy addition lemma}
Let $u : I^3\to X$ be a singular cube in a Hausdorff space $X$.
Then by restricting $u$ to the faces of $I^3$ and taking the
corresponding elements in $\brho^{\square}_2 (X)$, we obtain a
cube in $\brho^{\square} (X)$ which is commutative by the homotopy
addition lemma for $\brho^{\square} (X)$ (\cite{BHKP}, Proposition
5.5). Consequently, if $f : \brho^{\square} (X)\to \mathsf{D}$ is
a morphism of double groupoids with connections, any singular cube
in $X$ determines a commutative  {3-shell} in $\mathsf{D}$.

\section{The van Kampen/coequaliser Theorem}

The general setting of the van Kampen/coequaliser theorem is that
of a local-to global problem which can be explained as follows:

\begin{quote} Given an open covering $\mathcal{U}$ of $X$ and knowledge of
each ${\brho}^\square(U)$ for $U$ in $\mathcal{U}$, give a
determination of $\hdgb(X)$. \end{quote} Of course we need also to
know the values of $\hdgb$ on intersections $U \cap V$ and on the
inclusions from $U \cap V$ to $U$ and $V$.

We first note that that the functor $\hdgb$ on {the category
$\mathbf{Haus}$ of Hausdorff spaces } preserves coproducts
$\bigsqcup$, since these are just disjoint union in {Hausdorff
spaces} and in double groupoids. It is an advantage of the
groupoid approach that the coproduct of such objects is so simple
to describe.

Suppose we are given a cover $\cU$ of $X$. Then the homotopy
double groupoids in the following $\boldsymbol{\rho}$-{\em
sequence of the cover} are well-defined:

\begin{equation}
\bigsqcup _{(U,V)  \in \mathcal{U}\, ^2}
 \hdgb(U\cap V)\;
 \overset{a}{\underset{b}{\rightrightarrows}}  \bigsqcup _{U \in \cU}\hdgb(U)
 \labto{c} \hdgb(X ) .\label{coeq} \end{equation}
The morphisms      $a,b$ are determined by the inclusions
$$a_{UV}:U\cap V\rightarrow U, b_{UV}:U\cap V\rightarrow V$$ for
each $(U,V)  \in \mathcal{U}\,^2$ and $c$ is determined by the
inclusion $c_U:U\rightarrow X$ for each $U \in \cU$.

\begin{theorem}[van Kampen/coequaliser theorem]\label{twovkt} If  the
interiors of the sets of $\cU $ cover $X$, then in the above
$\boldsymbol{\rho}$-sequence of the cover, $c$ is the coequaliser
of $a,b$ in the category of double groupoids with connections.
\end{theorem}

A special case of this result is when $\cU$ has two elements. In
this case the coequaliser reduces  to a pushout.

The proof of the theorem is a direct verification of the universal
property for the coequaliser. So suppose $\mathsf{D}$ is a double
groupoid, and $$f:\bigsqcup _{U \in \cU}\hdgb(U) \to \mathsf{D}$$
is a morphism such that $fa=fb$. We require to construct uniquely
a morphism $F: \hdgb(X ) \to\mathsf{D}$ such that $Fc=f$. It is
convenient to write $f^U$ for the restriction of $f$ to
$\hdgb(U)$.

First, $F$ is uniquely defined on objects, since if $x \in X$ then
$x \in U$ for some $U \in \cU$ and so $F(x)=Fc(x)=f^U (x)$, and
the condition $fa=fb$ ensures this value is independent of $U$.

We next consider $F$ on an element $\langle u \rangle \in
\hdgb_1(X)$. By the Lebesgue covering lemma, we can write
$$ \< u \> = \< u_1\> + \cdots + \< u_n \> $$
where $u_i$ is a path in a set $U_i$ of the cover, and so
determines an element $\<u'_i\>$ in $\hdgb_1(U_i)$. The rule
$fa=fb$ implies that the elements $f^{U_i}\<u'_i\>$ are composable
in $\mathsf{D}$, and, since $F$ is a morphism, and $Fc=f$, their
sum must be $F\<u\>$.

A similar argument applies to an element $\<\alpha\> \in
\hdgb_2(X)$, where this time we choose a subdivision of $\alpha$
as a multiple composition $[\alpha_{ij}]$ with $\alpha_{ij}$ lying
in some $U_{ij}$. We must have $F\<\alpha\>=
[f^{{ij}}\<\alpha_{ij}'\>]$.

Now we must prove that $F$ can be well defined by choices of this
kind. Independence of the subdivision chosen is easily verified by
superimposing subdivisions. The hard part is to show the effect of
a homotopy.

We start in dimension 1.

Suppose $h : u \sim_T v$ is a thin homotopy of paths $u,v$. Then
there is a factorisation $$h: I ^2 \labto{\Phi} J \labto{p} X$$
where $J$ is a tree, and $\Phi$ is PWL on the boundary $\partial
I^2$ of $I^2$.

Now $p^{-1}(\cU)$ covers $J$. Choose a subdivision of $J$ so that
the open star of each edge of the subdivision is contained in a
set of $p^{-1}(\cU)$.

By a {\it grid subdivision of $I^n$} we mean a subdivision
determined by equally spaced hyperplanes of the form $x_i=$
constant, $i=1, \ldots, n$. Choose a grid subdivision of $I^2$
such that $\Phi$ maps each little square $s$ of the subdivision
into some open star $St(s)$ of an edge of $J$, and so $\Phi(s)
\subseteq p^{-1} U(s)$ for some $U(s) \in \cU$.

Assign extra vertices to $J$ as $\Phi$ of the vertices of the grid
subdivision.

Consider a non boundary edge $e$ of the grid subdivision, with
adjacent squares $s,s'$. Deform $\Phi$ on $e$, keeping its image
in $St(s) \cap St(s')$ so that $\Phi$ is PWL on $e$. Doing this
for each edge of a square $s$ gives a homotopy of the restriction
$\Phi\mid \partial s$ with image always in $St(s)$. The HEP allows
this homotopy  to be extended to a homotopy of $\Phi \mid s $
still with values in $St(s)$. These determine a homotopy $\Phi
\simeq \Phi'$ rel $\partial I^2$ such that if $h'= p\Phi'$, then
$h'\mid s$ is a thin square lying in $U(s) \in \cU$ for each $s$
of the grid subdivision. The class  $\< h'\mid s\> $ in $\hdgb_2
(U(s))$ is mapped by  $f^{U(s)}$ to a thin square in $\sD$. The
condition $fa=fb$ implies these squares are composable, and by the
property (T2) of thin squares,  their composite is a thin square.
By (T3) of \S 2, this composite is an identity, and so $F(\<u\>) =
F(\< v \> ) $. Thus  $F$ is well defined.

The proof that $F$ in this dimension is a morphism of groupoids
is immediate from the definition.

We apply a similar argument in dimension 2.

Suppose $W: \alpha \equiv_T \beta$ is a cubically thin homotopy.
Make a grid subdivision of $I^3$ into subcubes $c$ such that $W$
maps $c$ into a set $U(c) \in \cU$. Then $c$ determines a
commutative 3-shell in $\hdgb(U(c))$ which is mapped by $f^{U(c)}$
to a commutative 3-shell in $\sD$. Again by $fa=fb$ but now in the
next dimension, these commutative 3-shells are composable in $\sD$
to give a commutative 3-shell $C$. But the faces of $W$ not in
direction 1 are given to be thin homotopies. The argument as above
shows that the faces of $C$ not in direction 1 are thin squares in
$\sD$. By Lemma 2.6, $F(\<\alpha \>) = F(\< \beta \> ),$ and so
$F$ is well defined.

It is easy to check that $F$ is a morphism of double groupoids
with connections. This completes the proof of our main result,
apart from the proof of Theorem \ref{compcubes}.
\section{Proof of Theorem \ref{compcubes}}

\label{comproof}

The proof uses some 2-dimensional rewriting using connections of
the type used since the 1970s. However there are some tricky
points which we would like to emphasise.

We often have to rearrange some block subdivision of a multiple
composition. The general validity of this process is discussed in
\cite{DP1}, and its application to double categories with
connections in \cite{BM}. Here we point out the following.

Let $\alpha, \beta$ be squares in a double category such that
$$\gamma = \begin{bmatrix} \alpha& \beta\end{bmatrix} =  \alpha +_2
\beta$$ is defined. Suppose  further that $$ \alpha =
\begin{bmatrix} \alpha_1 \\ \alpha_2 \end{bmatrix} \quad \beta =
\begin{bmatrix} \beta_1 \\ \beta_2 \end{bmatrix}.$$ If $\alpha_1
+_2 \beta_1, \quad \alpha_2 +_2 \beta_2$ are defined, then we can
write $$ \gamma =\begin{bmatrix} \alpha_1& \beta_1 \\ \alpha_2 &
\beta_2
\end{bmatrix}.$$

However, if we rewrite $$\alpha =
\begin{bmatrix} \alpha_1' \\ \alpha_2' \end{bmatrix}, \quad \beta =
\begin{bmatrix} \beta_1' \\ \beta_2' \end{bmatrix}$$ then we
cannot write
$$ \gamma =\begin{bmatrix} \alpha_1'& \beta_1' \\ \alpha_2' &
\beta_2'
\end{bmatrix}$$ unless we are sure the compositions $\alpha_1'
+_2 \beta_1', \quad \alpha_2' +_2 \beta_2'$ are defined. Thus care
is needed in 2-dimensional rewriting.

We now proceed with the proof that all compositions of commutative
3-cubes are commutative. We first do the case of direction 2. This
requires the transport laws.

 Let $\alpha,\, \beta$ be  3-shells in $\mathsf K$ such that
 $\alpha^+_2 = \beta ^-_2$. We assume that $\alpha, \, \beta $ are
 commutative, so that the HCL holds for each. Then
 \begin{alignat*}{2}
\boldsymbol{\partial}^{\mathrm{odd}}(\alpha +_2 \beta) &=
\begin{bmatrix}
\; \tl & (\alpha +_2 \beta)^-_1 & \br \; \\[0.5ex]
 (\alpha +_2 \beta)^-_3 & (\alpha +_2 \beta)^+_2 & \hh
\end{bmatrix} \quad & & (\text{by definition})\\
&&&\\ &=\begin{bmatrix} \; \tl & \alpha^-_1  +_1 \beta^-_1 & \br
\; \\[0.5ex] \alpha^-_3  +_2 \beta^-_3 & \beta^+_2 & \hh
\end{bmatrix} && (\text{by definition of $+_2$ for 3-shells}) \\[0.5ex]
&=\begin{bmatrix} \; \tl & \hh & \alpha^-_1 & \br & \vv \\[0.5ex]
\vv & \tl & \beta^-_1 &  \hh&  \br \;
\\[0.5ex] \alpha^-_3 &  \beta^-_3 & \beta^+_2 & \hh &\hh
\end{bmatrix} && (\text{by the transport laws}) \\[1ex]
&=\begin{bmatrix} \; \tl & \hh & \alpha^-_1 & \br & \vv \\[0.5ex]
\alpha^-_3 & \hh & \beta^-_2 &  \hh&  \beta^+_3 \;
\\[0.5ex]  \vv & \tl & \beta^+_1 &  \hh &\br
\end{bmatrix} && \text{(by HCL for  $\beta$)}\\[0.5ex]
&=\begin{bmatrix} \; \tl & \hh & \alpha^-_1 & \br & \vv \\[0.5ex]
\alpha^-_3 & \hh & \alpha^+_2 &  \hh&  \beta^+_3 \;
\\[0.5ex]  \vv & \tl & \beta^+_1 &  \hh &\br
\end{bmatrix} && \text{(since $\beta^-_2=\alpha^+_2$ )}\\[0.5ex]
&=\begin{bmatrix} \; \hh & \hh & \alpha^-_2 & \alpha^+_3 & \beta^+_3 \\[0.5ex]
\tl & \hh & \alpha^+_1 &  \br&  \vv \;
\\[0.5ex]  \vv & \tl & \beta^+_1 &  \hh &\br
\end{bmatrix} && (\text{by HCL for $\alpha$}) \\
&= \begin{bmatrix} \hh & (\alpha +_2 \beta)^-_2 & (\alpha +_2
\beta)^+_3 \\[0.5ex]
\tl & (\alpha+_2 \beta) ^+_1 & \br \end{bmatrix} & & (\text{by
transport laws and  composition rules})\\[0.5ex]
&= \boldsymbol{\partial}^{\mathrm{even}}(\alpha +_2 \beta) &&
(\text{as required.})
\end{alignat*}

We now consider the case of $+_3$ and for this it turns out to be
convenient to use (HCL$'$). We write the left hand side of this as
$ \boldsymbol{\partial}'^{\mathrm{odd}}$ and the right hand side
as $ \boldsymbol{\partial}'^{\mathrm{even}}$. Suppose then $
\alpha +_3 \beta$ is defined. Then
\begin{alignat*}{2}
\boldsymbol{\partial}'^{\mathrm{odd}}(\alpha +_3 \beta) &=
\begin{bmatrix} \tl & (\alpha +_3 \beta)^-_1 \\
(\alpha +_3 \beta)^-_3 &(\alpha +_3 \beta)^+_2\\
\br & \vv
\end{bmatrix} && \\
&= \begin{bmatrix} \tl & \alpha^-_1 +_2  \beta^-_1 \\
\alpha^-_3 &\alpha^+_2  +_2 \beta^+_2\\
\br & \vv
\end{bmatrix}
&&\\ &= \begin{bmatrix} \tl & \alpha^-_1 &   \beta^-_1 \\
\alpha^-_3 &\alpha^+_2 &  \beta^+_2\\
\br & \vv & \vv
\end{bmatrix} & \\
&= \begin{bmatrix} \vv & \tl &   \beta^-_1 \\
\alpha^-_2 &\alpha^+_3 &  \beta^+_2\\
\alpha^+_1 &{\br} & \vv
\end{bmatrix} && \text{(by HCL}') \\
&=\begin{bmatrix} \vv & \tl &   \beta^-_1 \\
\alpha^-_2 &\beta^-_3 &  \beta^+_2\\
\alpha^+_1  & \br  & \vv
\end{bmatrix}&& (\text{as } \alpha^+_3 = \beta^-_3 ) \\
&=\begin{bmatrix} \vv &\vv &  \tl &   \\
\alpha^-_2 &\beta^-_2 &  \beta^+_3\\
\alpha^+_1  & \beta^+_1  & {\br}
\end{bmatrix} && (\text{by HCL}') \\
&=\begin{bmatrix} \vv &  \tl    \\
(\alpha +_3 \beta)^-_2 &  (\alpha +_3\beta)^+_3\\
(\alpha  +_3 \beta) ^+_1  & {\br}
\end{bmatrix} && \\&= \partial'^{\mathrm{even}} (\alpha +_3
\beta)&& \mbox{(as required.)}
\end{alignat*}

The proof for $+_1$ is similar to the last one, but using (HCL).
We leave it to the reader.
 \mbox{} \hfill $ \blacksquare $

To see the complications of these ideas in higher dimensions, see
\cite{ABS,Higgins} and also p.361-362 of \cite{Gaucher}, which
deals with the 4-cube.

Note also that these calculations, when applied to
$\boldsymbol{\rho}^\square(X)$, imply the existence, even a
construction,  of certain homotopies which would otherwise be
difficult to find. The possibility of calculating with such
homotopies is indeed one of the aims of this theory.
\section{Future directions}
Our intention for future papers in this series is to investigate:
\begin{enumerate}[A)]
\item Does $\hdgb(X)$ capture the weak homotopy 2-type of $X$?
This seems likely in view of the facts on $\hdgb(X)$ given in the
Introduction, and since crossed modules of groupoids, and hence
also double groupoids with connections, capture all weak homotopy
2-types (see, for example, \cite{class}).
 \item How does the functor $\hdgb$ behave on homotopies? This
 question requires for its answer notions of
 tensor product and homotopies for double groupoids with
 connections. These have been developed in the corresponding
 $\omega$-groupoid case in \cite{tens}, and in the
 $\omega$-category case in \cite{ABS}.

 \item It is also expected that the above developments will allow
 an enrichment of the category of Hausdorff spaces over a monoidal
 closed category of double groupoids with connection, analogous to
 the enrichment of the category of filtered spaces over the
 category of crossed complexes, given in \cite{class}. Such an
 enrichment should be useful for calculations with some multiple
 compositions of homotopies.

 \item Are there smooth analogues of $\hdgb(M)$ in the case when
 $M$ is a smooth manifold? This possibility is suggested by the
 methods of \cite{mackaay}.

\end{enumerate}

\end{document}